\documentclass{article}
\usepackage{iclr2027_conference,times}

\usepackage[utf8]{inputenc} 
\usepackage[T1]{fontenc}    
\usepackage{hyperref}       
\usepackage{url}            
\usepackage{booktabs}       
\usepackage{amsfonts}       
\usepackage{nicefrac}       
\usepackage{microtype}      
\usepackage{xcolor}         
\usepackage{amsthm}
\usepackage{amsmath}
\usepackage{amssymb}
\usepackage{mathrsfs}
\usepackage{multirow}
\usepackage{graphicx}
\usepackage{subcaption}
\usepackage{caption}
\usepackage{diagbox}
\usepackage{array}
\usepackage{makecell}
\newcolumntype{C}[1]{>{\centering\arraybackslash}m{#1}}
\usepackage[ruled]{algorithm2e}

\newtheorem{theorem}{Theorem}
\newtheorem{lemma}{Lemma}

\newtheorem{assumption}{Assumption}
\usepackage{amsmath,amssymb,amsfonts,amsthm,mathrsfs}
\usepackage{booktabs,array,longtable}
\makeatletter
\@ifundefined{theorem}{\newtheorem{theorem}{Theorem}}{}
\@ifundefined{lemma}{}{}
\@ifundefined{proposition}{}{}
\@ifundefined{corollary}{}{}
\@ifundefined{assumption}{}{}
\@ifundefined{restatedtheorem}{}{}
\makeatother

\setkeys{Gin}{draft=false}

\title{From What Is Used to What Is Missing:\\ Rethinking Dynamic Rank in LoRA}

\author{{Peng Xue}\textsuperscript{1}\quad{Guodong Fan}\textsuperscript{2}\quad
{Guang-Yong Chen}\textsuperscript{1}\quad
{Xing Chen}\textsuperscript{1}\\
\normalfont\textsuperscript{1} Fuzhou University\\
\normalfont\textsuperscript{2} Shandong Technology and Business University
}

\hypersetup{
  pdftitle={From What Is Used to What Is Missing: Rethinking Dynamic Rank in LoRA},
  pdfauthor={Peng Xue, Guodong Fan, Guang-Yong Chen, Xing Chen},
  pdfsubject={Under review as a conference paper at ICLR 2027}
}

\iclrfinalcopy

\begin{document}

\maketitle
\lhead{Under review as a conference paper at ICLR 2027}

\begin{abstract}
Parameter-efficient fine-tuning with LoRA is governed by rank, which controls parameter budget and locally accessible task updates. Existing dynamic-rank methods mainly focus on capacity allocation, leaving unresolved whether the current adaptation space is sufficient for the task. We characterize the first-order information boundary of standard factor gradients, showing that their unobservable subspace coincides with the normal space of the low-rank manifold. Thus, factor gradients alone cannot distinguish under-optimization within the current space from task-relevant directions beyond its first-order reach. This motivates a geometric view of rank evolution, where tangent information optimizes existing degrees of freedom and the normal spectrum quantifies the value of new ones. Under a local quadratic model, we show that leading normal singular directions yield the locally optimal rank expansion, with spectral energy determining its gain. Building on this characterization, we propose GEAR, which turns the tangent--normal decomposition into rank evolution through tangent-space optimization, matrix-free normal-spectrum estimation for rank expansion, and task-induced geometric initialization. Across language and vision benchmarks, GEAR achieves the best average performance in all evaluated settings with substantially fewer trainable parameters. On LLaMA3-8B, it raises the average commonsense score from 86.61 to 88.45 while reducing trainable parameters from 42.0M to 14.4M.
\end{abstract}

\section{Introduction}
Parameter-efficient fine-tuning (PEFT) has become a prevailing approach to adapting large pretrained models by updating only a small fraction of task-specific parameters \citep{ding2023peft}. Among PEFT methods, LoRA represents the weight increment as a low-rank matrix, \(\Delta W=BA\), where \(A\in\mathbb{R}^{r\times m}\) and \(B\in\mathbb{R}^{n\times r}\), reducing the number of trainable parameters from $nm$ to $r(n+m)$ \citep{hu2022lora}. While \(r\) therefore directly controls the adaptation budget, it also determines the local first-order update space accessible under the low-rank constraint \citep{vandereycken2013lowrank}. Rank thus governs both adaptation capacity and which task-induced directions can be expressed. Since these requirements vary across layers and throughout training, a uniform fixed rank is unlikely to match the evolving structure of task demand, motivating dynamic-rank adaptation.

Most dynamic-rank methods still treat rank primarily as a capacity-allocation variable. AdaLoRA, for example, reallocates a prescribed rank budget across weight matrices according to the importance of individual singular components \citep{zhang2023adalora}, while DyLoRA trains nested low-rank structures to support inference at multiple truncation ranks \citep{valipour2023dylora}. More recent approaches permit more flexible rank evolution: FlexLoRA uses the spectral entropy of the current low-rank update to determine when capacity should be pruned or expanded \citep{liu2026flexlora}, whereas GeoLoRA combines gradient-driven basis augmentation with spectral truncation within a dynamical low-rank integration framework \citep{schotthofer2025geolora}. Although these methods substantially improve the flexibility of rank adaptation, their rank decisions are primarily derived from the capacity or structure already represented by the current adapter. They therefore characterize what the current low-rank space uses, but provide no direct criterion for determining whether task-relevant descent directions remain missing outside that space. This leaves a more fundamental question unresolved:

\textit{Is the current adaptation space locally sufficient to express task-relevant descent directions, or must its rank expand to make such directions first-order accessible?}

Answering this question requires distinguishing insufficient optimization within the current low-rank space from insufficiency of the space itself. Limited progress may arise for two fundamentally different reasons: the current rank may be sufficient but under-optimized, or task-relevant descent directions may lie outside its first-order reachable space. Only the latter justifies increasing rank. Existing rank-allocation criteria do not explicitly separate these two cases, making it difficult to determine whether additional capacity is genuinely required. This distinction motivates our central view that dynamic rank should be driven not only by what the current adaptation space already represents, but also by what task-relevant descent information remains inaccessible to it.

\begin{figure}[htbp]
\hypertarget{fig:advance}{\mbox{}}
\centering
{\includegraphics[width=\linewidth]{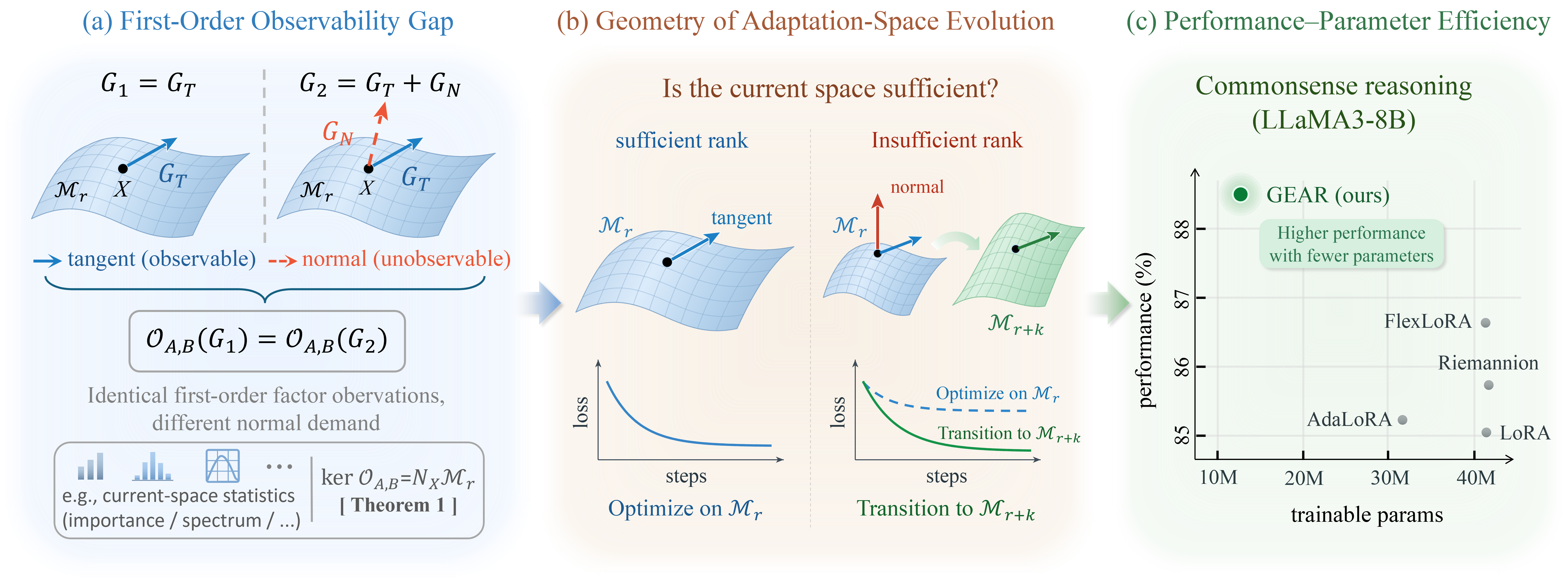}}
\caption{
Motivation and overview of GEAR.
(a) standard factor-gradient observations are insensitive to normal components, so \textbf{identical first-order observations can correspond to different rank requirements}.
(b) GEAR distinguishes optimization within the current manifold from expansion toward task-relevant directions beyond the current first-order reachable space.
(c) this geometry-guided rank evolution yields a stronger performance--parameter trade-off on commonsense reasoning.
}
\end{figure}

To address this limitation, we propose GEAR (Geometry-guided Evolution of Adaptive Rank), which formulates dynamic-rank adaptation as the evolution of a task-dependent low-rank space. Rather than inferring rank adequacy from existing factor statistics, we characterize the first-order information available under the fixed-rank parameterization. For \(X=BA\), we show that the standard factor-gradient map has an observability gap whose null space coincides with the normal space \(N_X\mathcal M_r\) of the fixed-rank manifold \(\mathcal M_r=\{X\in\mathbb R^{n\times m}:\operatorname{rank}(X)=r\}\). This establishes a structural information limitation in which standard factor gradients reveal the tangent component relevant to optimizing the current rank, while leaving the normal component unobservable. Consequently, they cannot distinguish under-optimization within the current adaptation space from task-relevant directions requiring rank expansion to become first-order accessible. For a prescribed rank increase, we further show that leading normal modes determine the locally optimal expansion, while their spectral energy quantifies its gain. Building on this characterization, we propose GEAR, which turns the tangent--normal decomposition into a practical rank-evolution mechanism through tangent-space optimization, matrix-free normal-spectrum estimation for rank expansion, and task-induced geometric initialization. By expanding rank only when the estimated normal spectrum reveals valuable task-relevant directions, GEAR improves adaptation effectiveness and parameter efficiency. Across multiple benchmarks, it consistently outperforms competitive PEFT baselines while using fewer trainable parameters. Our contributions are as follows:
\begin{itemize}
    \item \textbf{First-order observability barrier for rank adaptation.} To our knowledge, we are the first to characterize the information available from standard fixed-rank factor gradients and prove that their unobservable subspace is the normal space of the low-rank manifold, thereby separating current-space optimization from rank-expanding directions.
    \item \textbf{Geometry-guided rank evolution.} We introduce GEAR, which unifies optimization over the current fixed-rank manifold and rank expansion using tangent information to optimize the adaptation space and the normal spectrum to identify and value new task-relevant directions.
    \item \textbf{Excellent performance–parameter trade-off.} Across language and vision benchmarks, GEAR achieves the best average performance in every evaluated setting while using substantially fewer trainable parameters than strong PEFT baselines.
\end{itemize}

\section{Related work}
\label{sec:2}
\paragraph{LoRA and Adaptive Rank Adaptation.} 
LoRA parameterizes task-specific updates through low-rank factors, reducing trainable parameters and storage while preserving frozen pretrained weights \citep{hu2022lora}. Subsequent work has improved low-rank adaptation through refined parameterizations and initialization strategies. DoRA decouples weight magnitude and direction to enhance adaptation flexibility \citep{liu2024dora}, while PiSSA initializes adapters from principal singular subspaces of the pretrained weights \citep{meng2024pissa} and LoRA-GA uses task gradients to identify more effective initial low-rank directions \citep{wang2024loraga}. Despite these advances, most methods still operate with a fixed rank, which may poorly match the heterogeneous adaptation demands across layers and training stages. This limitation has motivated adaptive-rank LoRA methods. AdaLoRA dynamically reallocates a prescribed rank budget across modules according to singular-component importance \citep{zhang2023adalora}; DyLoRA trains nested low-rank structures to support multiple effective ranks \citep{valipour2023dylora}; and SoRA learns effective rank through sparse gating \citep{ding2023sora}. More recent approaches allow rank to both contract and expand: FlexLoRA adjusts rank according to the spectral structure of the current adapter \citep{liu2026flexlora}, whereas GeoLoRA combines gradient-driven basis augmentation with spectral truncation under a dynamical low-rank formulation \citep{schotthofer2025geolora}. GeLoRA further estimates the intrinsic dimensionality of hidden representations to derive adaptive layer-wise ranks \citep{eddib2025gelora}. Collectively, these methods have advanced LoRA from fixed-rank adaptation toward increasingly flexible, task-dependent capacity allocation. However, existing adaptive-rank methods primarily decide how to allocate or reorganize represented low-rank capacity, whereas GEAR asks a different question: whether task-relevant descent information remains outside the current first-order reachable space and therefore warrants new rank dimensions.

\paragraph{Geometric Optimization for LoRA}
Fixed-rank matrices admit a natural Riemannian manifold structure, providing a principled way to optimize low-rank states independently of a particular factorization \citep{vandereycken2013lowrank}. Early geometric approaches to LoRA exploit this structure through Riemannian metrics and preconditioning; for example, Riemannian Preconditioned LoRA derives low-dimensional preconditioners from the geometry of the low-rank matrix space to improve optimization stability and convergence \citep{zhang2024riemannian}. GeoLoRA further introduces dynamical low-rank integration into PEFT, combining gradient-driven basis augmentation with spectral truncation for structure-preserving and rank-adaptive updates \citep{schotthofer2025geolora}. Related work has also explored explicit manifold constraints on LoRA factors. Park et al. optimize the \(B\) factor on the Stiefel manifold to maintain orthogonality and effective rank \citep{park2025riemannian}, while Manifold-LoRA develops a retraction-free Stiefel-manifold formulation for LoRA fine-tuning \citep{zhang2026manifoldlora}. More recently, Riemannion directly optimizes the adaptation matrix on the fixed-rank manifold, removing factorization ambiguity and enabling parametrization-independent geometric updates \citep{bogachev2026riemannion}.

These methods establish geometry as an effective tool for improving low-rank optimization, primarily by refining the metric, constraining the factor geometry, or evolving the current low-rank state. GEAR instead focuses on the information required to determine when the current adaptation space itself should expand. From a first-order observability perspective, it relates the normal component unobservable from standard factor gradients to task-relevant directions beyond the current first-order reachable space, thereby connecting manifold optimization with rank expansion.

\section{Method}
\label{sec:3}
GEAR builds on a first-order observability gap in fixed-rank adaptation: standard factor gradients characterize optimization within the current space, but not the normal task information that governs its expansion. We first formalize this boundary and then develop a geometry-guided adaptation framework that optimizes the low-rank space while identifying when and how it should expand.

\begin{algorithm}[t]
\caption{GEAR: Geometry-Guided Evolution of Adaptive Rank}
\label{alg:gear}

\KwIn{
Pretrained weights $W$;
initial rank $r_0$;
rank upper bound $r_{\max}$;
audit set $\mathcal T_{\mathrm{audit}}$;
budget $\mathcal B_{\mathrm{target}}$
}

\BlankLine
$(U,C,V) \leftarrow \mathrm{TGI}(W,r_0)$\;
$X \leftarrow UC,\quad W \leftarrow W-sX,\quad r \leftarrow r_0$\;

\For{$t=1,\ldots,T$}{

    $G_T \leftarrow
    U(U^\top G)+(I-UU^\top)GVV^\top$
    \tcp*[r]{Eq.~\eqref{eq:tangent_grad}}

    $(X, U, C, V) \leftarrow
    \operatorname{Retr}_{\mathcal M_r}(X-\eta_t G_T)$
    \tcp*[r]{Eq.~\eqref{eq:manifold_update}}

    \If{$t\in\mathcal T_{\mathrm{audit}}$ \textbf{and} $r<r_{\max}$}{

        $(\widehat U,\widehat\Sigma,\widehat V)
        \leftarrow
        \mathrm{NormalSketch}(X)$
        \tcp*[r]{Eqs.~\eqref{eq:normal_probe}--\eqref{eq:normal_sketch}}

        $\mathcal S
        \leftarrow
        \mathrm{Select}(
        \widehat U,\widehat\Sigma,\widehat V;
        \mathcal B_{\mathrm{target}})$
        \tcp*[r]{Eq.~\eqref{eq:rank_score}}

        $X_{\mathrm{cand}}
        \leftarrow
        \mathrm{Propose}(X,\mathcal S)$
        \tcp*[r]{Eq.~\eqref{eq:candidate_expansion}}

        \If{$\mathrm{Validate}(X_{\mathrm{cand}},X)$}{

            $(U,C,r)
            \leftarrow
            \mathrm{Expand}(U,C,\mathcal S),\quad X\leftarrow UC$
            \tcp*[r]{Eq.~\eqref{eq:rank_expansion}}
        }
    }
}
\KwOut{$X$}
\end{algorithm}

\subsection{First-Order Information Boundary of Standard Factor Gradients}
Consider a pretrained weight matrix \(W\in\mathbb{R}^{n\times m}\) adapted as
$$ W'=W+sX,\qquad X=BA, $$
where \(s>0\) is a fixed scaling factor. Let \(f(X)=\mathcal L(W+sX)\) and \(G=\nabla_X f(X)\). For a rank-\(r\) adaptation matrix \(X=BA\), we analyze \(X\) on the fixed-rank manifold \(\mathcal M_r\), whose tangent space characterizes the first-order directions accessible without increasing rank \citep{vandereycken2013lowrank}.

Let \(X=U\Sigma V^\top\) be the compact SVD of \(X\), and define \(P_U=UU^\top\) and \(P_V=VV^\top\). The full gradient admits the orthogonal decomposition
$$ G=G_T+G_N, $$
where $ G_T=P_UG+(I-P_U)GP_V\in T_X\mathcal M_r$ and $G_N=(I-P_U)G(I-P_V)\in N_X\mathcal M_r$. Under the factorized parameterization \(X=BA\), the corresponding gradients are
$$ \nabla_B f=GA^\top,\qquad \nabla_A f=B^\top G. $$
We therefore define the factor-gradient observation operator
$$ \mathcal O_{A,B}(G) = \bigl(GA^\top,B^\top G\bigr), $$
which captures the first-order information about \(G\) available from standard factor backpropagation.
\begin{theorem}[Factor-Gradient Observability Barrier]
\label{thm:1}
Let \(A\in\mathbb R^{r\times m}\) and \(B\in\mathbb R^{n\times r}\) be full rank, with \(X=BA\in\mathcal M_r\). Then
$$ \ker\mathcal O_{A,B} = N_X\mathcal M_r. $$
Equivalently, for any \(G_1\) and \(G_2\),
$$ \mathcal O_{A,B}(G_1) = \mathcal O_{A,B}(G_2) $$
if and only if $ (G_1)_T=(G_2)_T.$
\end{theorem}

The theorem establishes a precise first-order observability boundary, showing that standard factor gradients preserve the tangent component of the task gradient but cannot distinguish gradients that differ only in their normal components. A complete proof is provided in Appendix~\ref{app:proof_theorem1}. Consequently, the same factor-gradient observation may correspond either to a gradient fully represented by the current adaptation space or to one containing additional task-relevant directions beyond its first-order reach. This ambiguity arises from the information available under the current fixed-rank parameterization, rather than from factorization itself.

This result distinguishes optimization within the current adaptation space from the need for rank expansion. The tangent component represents directions that are first-order accessible at rank \(r\), whereas the normal component captures task-relevant directions beyond the current space. Accordingly, tangent information is used to optimize the existing degrees of freedom, while normal information provides the relevant signal for assessing the value of rank expansion. The subsequent sections develop these two mechanisms, respectively.

\subsection{Fixed-Rank Manifold Optimization}
\label{sec:3.2}
For a fixed rank \(r\), GEAR performs optimization on the manifold \(\mathcal M_r\) using the tangent component of the gradient. We represent the adaptation matrix as
$$
X=UC,\qquad U^\top U=I_r,
$$
where \(U\) is an orthonormal basis for the column space of \(X\), and \(C\) contains the corresponding low-rank coefficients. Let \(V\) be an orthonormal basis for the column space of \(C^\top\). The tangent gradient is then
\begin{equation}
G_T = U(U^\top G) + (I-UU^\top)GVV^\top.
\label{eq:tangent_grad}
\end{equation}
Equation \eqref{eq:tangent_grad} shows that \(G_T\) can be computed from the two gradient contractions \(U^\top G\) and \(GV\), avoiding explicit construction of the full gradient \(G\).

To obtain these quantities directly through backpropagation, GEAR uses the forward-equivalent representation
$$
B_{\mathrm{ext}} = [U,\ 0], \qquad A_{\mathrm{ext}} = \begin{bmatrix} C\\ V^\top \end{bmatrix},
$$
which satisfies \(B_{\mathrm{ext}}A_{\mathrm{ext}}=X\). Although this representation leaves the forward mapping unchanged, its backward derivatives expose the required gradient contractions:
$$
\nabla_C f=U^\top G, \qquad \nabla_{B_{\mathrm{ext}}^{(2)}}f=GV.
$$
Thus, the tangent gradient can be recovered without explicitly materializing \(G\). GEAR then performs a tangent update followed by retraction onto the fixed-rank manifold,
\begin{equation}
X^{-}=X-\eta G_T,
\qquad
X_{t+1}=\operatorname{Retr}_{\mathcal M_r}(X^{-}).
\label{eq:manifold_update}
\end{equation}
Since \(X^{-}\) has rank at most $2r$, the retraction is performed on its compact factors rather than on the full matrix. Specifically, thin QR factorizations reduce the update to a small core matrix, whose truncated SVD yields the retracted rank-\(r\) state together with the updated left and right subspace bases. Full details are provided in Appendix~\ref{app:detail_retraction}.

\subsection{Geometry-Guided Rank Evolution}
The preceding section optimizes the current adaptation state on \(\mathcal M_r\). Rank expansion, by contrast, is governed by the normal component \(G_N\), which captures task-relevant directions beyond the current first-order reachable space. We next quantify the local value of these directions for rank expansion.
Let
$$
G_N=\sum_{i\ge 1}\sigma_i u_i v_i^\top, \qquad \sigma_1\ge\sigma_2\ge\cdots,
$$
and denote its best rank-$k$ approximation by $(G_N)_k$.
\begin{theorem}[Local Optimality of Rank Expansion]
\label{thm:2}
Let $X\in\mathcal M_r$, and let $\mathcal C_{r,k}(X)$ denote the tangent cone at X to the bounded-rank set $\mathcal V_{r+k} = \{Y\in\mathbb R^{n\times m}:\operatorname{rank}(Y)\le r+k\}$. For the local quadratic model $q_\eta(H) = \langle G,H\rangle_F + \frac{1}{2\eta}\|H\|_F^2$, an optimal direction over $\mathcal C_{r,k}(X)$ is
$$
H_k^\star = -\eta\left[G_T+(G_N)_k\right],
$$
where $(G_N)_k$ is a best rank-$k$ approximation of $G_N$. Moreover, the optimal model decrease enabled by the additional $k$ rank dimensions, relative to optimization over $T_X\mathcal M_r$, is
$$
\Delta_k := \min_{H\in T_X\mathcal M_r}q_\eta(H) - \min_{H\in\mathcal C_{r,k}(X)}q_\eta(H) = \frac{\eta}{2}\sum_{i=1}^{k}\sigma_i^2(G_N).
$$
\end{theorem}
A complete proof is provided in Appendix~\ref{app:proof_theorem2}. The result shows that the leading singular directions of the normal gradient constitute the locally optimal directions for rank expansion, while their spectral energy exactly quantifies the additional decrease attainable under the local quadratic model.

However, constructing \(G_N=(I-UU^\top)G(I-VV^\top)\) requires the full gradient \(G\). GEAR instead estimates its dominant spectrum through randomized probing \citep{halko2011randomized}. Given \(\Omega\in\mathbb R^{m\times h}\), with \(h\ll\min\{n,m\}\), we project it onto the orthogonal complement of the right subspace, \(\Omega_\perp=(I-VV^\top)\Omega\). We introduce a zero-initialized auxiliary variable \(P\) through \(\Delta X_{\mathrm{probe}}=P\Omega_\perp^\top\). Since \(P=0\), this leaves the forward mapping unchanged, while backpropagation yields \(\nabla_P f=G\Omega_\perp\). Projecting onto the orthogonal complement of the left subspace gives
\begin{equation}
Y=(I-UU^\top)\nabla_P f = G_N\Omega.
\label{eq:normal_probe}
\end{equation}
Let $Q=\operatorname{orth}(Y)$. If $G_N=U_N\Sigma_NV_N^\top$, then $Y=U_N\Sigma_N(V_N^\top\Omega)$, so singular directions associated with larger singular values contribute more strongly to $Y$. Hence, $\operatorname{span}(Q)$ captures the dominant left singular subspace of $G_N$. Applying the same zero-output probe yields \(Q^\top G\), from which we compute
\begin{equation}
Z_N
=
Q^\top G(I-VV^\top)
=
Q^\top G_N
=
\widetilde U\widehat\Sigma\widehat V^\top .
\label{eq:normal_sketch}
\end{equation}
The leading normal components are therefore approximated by $\widehat u_i=Q\widetilde u_i$, $\widehat v_i=\widehat V_{:,i}$, and $\widehat\sigma_i=\widehat\Sigma_{ii}$, without explicitly materializing $G_N$. The resulting approximation has a controlled effect on the local value of rank expansion.
\begin{theorem}[Stability under Inexact Normal Estimation]
\label{thm:3}
Let $\widehat G_N\in N_X\mathcal M_r$ satisfy $\|\widehat G_N-G_N\|_F\le\varepsilon$, and let $\tau_k=\|G_N-(G_N)_k\|_F$. Then the rank-k expansion based on $(\widehat G_N)_k$ satisfies
$$
0\le q_\eta\!\left(-\eta(\widehat G_N)_k\right) - q_\eta\!\left(-\eta(G_N)_k\right) \le 2\eta\varepsilon(\tau_k+\varepsilon).
$$
\end{theorem}
A complete proof is provided in Appendix~\ref{app:proof_theorem3}. GEAR then uses the estimated normal spectrum to allocate rank under a global parameter budget. Motivated by Theorem~\ref{thm:2}, we score each candidate normal direction by its estimated local gain per added parameter and impose a global rank budget:
\begin{equation}
\begin{aligned}
a_{\ell,i}
=
\frac{\widehat\sigma_{\ell,i}^{\,2}}
{n_\ell+m_\ell},
\qquad
\mathcal B_{\mathrm{target}}
=
\bar r_{\mathrm{budget}}
\sum_{\ell=1}^{L}
\left(n_\ell+m_\ell\right).
\end{aligned}
\label{eq:rank_score}
\end{equation}
Here, $r_\ell$ denotes the current rank of module $\ell$, $n_\ell+m_\ell$ is the cost of one additional rank dimension in module $\ell$, and $\bar r_{\mathrm{budget}}$ denotes the weighted-average rank budget. Candidates are selected according to $a_{\ell,i}$ subject to the global budget and the per-module rank limit. For \(k_{\ell}\) selected directions, let \(\widehat U_{\ell,k_{\ell}}=[\widehat u_{\ell,1},\ldots,\widehat u_{\ell,k_{\ell}}]\), \(\widehat V_{\ell,k_{\ell}}=[\widehat v_{\ell,1},\ldots,\widehat v_{\ell,k_{\ell}}]\), and \(\widehat\Sigma_{\ell,k_{\ell}}=\operatorname{diag}(\widehat\sigma_{\ell,1},\ldots,\widehat\sigma_{\ell,k_{\ell}})\). They define the candidate expansion
\begin{equation}
X_{\ell,\mathrm{cand}}
=
X_\ell-\eta_b
\widehat U_{\ell,k_{\ell}}
\widehat\Sigma_{\ell,k_{\ell}}
\widehat V_{\ell,k_{\ell}}^\top .
\label{eq:candidate_expansion}
\end{equation}
We validate \(X_{\ell,\mathrm{cand}}\) on a small set of separate training microbatches shared across candidate variants and accept the expansion only if its mean loss improves by more than a relative margin \(\rho\); full details are provided in Appendix~\ref{app:detail_validate}.

For an accepted expansion, $\widehat U_{\ell,k_{\ell}}\perp U_\ell$ and $\widehat V_{\ell,k_{\ell}}\perp V_\ell$, so the low-rank state is updated directly as
\begin{equation}
U_\ell^{+}=[U_\ell,\widehat U_{\ell,k_{\ell}}],
\qquad
C_\ell^{+}=
\begin{bmatrix}
C_\ell\\
-\eta_b\widehat\Sigma_{\ell,k_{\ell}}\widehat V_{\ell,k_{\ell}}^\top
\end{bmatrix},
\qquad
X_\ell^{+}=U_\ell^{+}C_\ell^{+}.
\label{eq:rank_expansion}
\end{equation}
The accepted expansion moves the adaptation state from \(\mathcal M_{r_\ell}\) to \(\mathcal M_{r_\ell+k_\ell}\). Optimizer states associated with existing components are inherited, while those of newly introduced directions are initialized to zero. Subsequent updates then follow the fixed-rank manifold optimization in Section~\ref{sec:3.2}.

\subsection{Task-Induced Geometric Initialization (TGI)}

Rank inadequacy is defined relative to the current adaptation space. To avoid conflating random subspace mismatch with genuine missing directions, GEAR initializes this space from task-gradient structure estimated on a small calibration set.
Let $G_0=\nabla_X f(0)$. In a left-probing stage, set
$X=P\Omega^\top$ with $P=0$ and fixed
$\Omega\in\mathbb R^{m\times h}$. This preserves the pretrained
mapping, while backpropagation yields
$Z_0=\left.\nabla_P f(P\Omega^\top)\right|_{P=0}
=G_0\Omega\in\mathbb R^{n\times h}$ without forming $G_0$.

To express the normalization for both stages uniformly,
write $Z_0\in\mathbb R^{d\times h}$, with $d=n$ for the left probe
and $d=m$ for the right probe.
Set $M=Z_0^{\odot2}$ and $\mu=(dh)^{-1}\sum_{i,a}M_{ia}$.
To prevent large gradient scales from dominating subspace estimation,
we normalize the sketch along both dimensions:
\[
\begin{aligned}
d_i^{\mathrm{row}}
  =\frac{1}{h}\sum_a M_{ia}+\lambda\mu,\quad
d_a^{\mathrm{col}}
  =\frac{1}{d}\sum_i M_{ia}+\lambda\mu,\quad
\widetilde Z_0=D_{\mathrm{row}}^{-1/2}Z_0D_{\mathrm{col}}^{-1/2},
\end{aligned}
\]
where $D_{\mathrm{row}}=\operatorname{diag}(d_i^{\mathrm{row}})$
and $D_{\mathrm{col}}=\operatorname{diag}(d_a^{\mathrm{col}})$. A thin QR factorization of the left-stage $\widetilde Z_0$ yields
$Q_L\in\mathbb R^{n\times h}$. Holding $Q_L$ fixed, a
zero-initialized right probe $X=Q_LR$ gives the transposed sketch
$Z_0=\left(\left.\nabla_R f(Q_LR)\right|_{R=0}\right)^\top
=G_0^\top Q_L\in\mathbb R^{m\times h}$.
Non-final right stages use QR to update the next left-probing basis.
TGI performs $p+1$ left--right rounds, totaling $2(p+1)$ stages.
Additional rounds can refine the dominant-subspace estimate by
probing with successively updated bases; we use $p=1$ in our experiments.

At the final right stage, the normalized sketch
$\widetilde Z_0\in\mathbb R^{m\times h}$ gives
$\bar Z_0=\widetilde Z_0^\top=\widehat U\Sigma V^\top$.
Retaining the leading $r_0\le h$ components, we initialize
\begin{equation}
U_0=Q_L\widehat U_{:,1:r_0},\qquad
C_0=\operatorname{diag}(\gamma_1,\ldots,\gamma_{r_0})
V_{:,1:r_0}^{\top},\qquad
X_0=U_0C_0.
\label{eq:tgi_init}
\end{equation}
Here $\gamma_i$ is the clipped ratio of a retained singular value
to the mean of the retained values. Since $X_0$ is generally nonzero, we compensate its initial
contribution in the frozen weight:
\[
\widetilde W=W-sX_0,\qquad
\widetilde W+sX_0=W.
\]
Thus, training starts from a task-aligned adaptation space while
preserving the pretrained forward mapping at initialization.

\section{Experimental results}
\label{sec:exp}
We evaluate GEAR across commonsense reasoning, natural language understanding, and visual recognition benchmarks, considering both predictive performance and parameter efficiency. We compare against representative fixed-rank and adaptive-rank PEFT methods, including LoRA \citep{hu2022lora}, AdaLoRA \citep{zhang2023adalora}, Riemannion \citep{bogachev2026riemannion}, GeoLoRA \citep{schotthofer2025geolora}, and FlexLoRA \citep{liu2026flexlora}. All methods use consistent experimental settings for fair comparison. Ablation studies are provided in Appendix~\ref{app:ablation}.

\subsection{Rank Audit Schedule}
To ensure stable rank adaptation, GEAR performs rank-growth audits periodically during early-to-middle training and freezes the rank thereafter. Let \(t_{\mathrm{start}}\) denote the first audit step, \(\Delta t\) the audit interval, and \(t_{\mathrm{final}}\) the step after which rank expansion stops. The audit steps are defined as
$$ \mathcal{T}_{\mathrm{audit}} = \left\{ t_{\mathrm{start}}+k\Delta t \;\middle|\; k\in\mathbb Z_{\ge 0}, \; t_{\mathrm{start}}+k\Delta t\le t_{\mathrm{final}} \right\}. $$
At each audit, GEAR estimates task-relevant directions beyond the current adaptation space and evaluates them as candidates for rank expansion. After \(t_{\mathrm{final}}\), the rank remains fixed and training proceeds solely with manifold optimization at the resulting rank.

\subsection{Commonsense Reasoning}
We evaluate GEAR on eight commonsense reasoning benchmarks, including BoolQ, PIQA, SIQA, and HellaSwag \citep{clark2019boolq,bisk2020piqa,sap2019socialiqa,zellers2019hellaswag}, together with WinoGrande, ARC-e/ARC-c, and OBQA \citep{sakaguchi2020winogrande,clark2018arc,mihaylov2018openbookqa}, using LLaMA3-8B \citep{grattafiori2024llama3} and Qwen2.5-7B \citep{qwen2025qwen25} as backbone models. The reference rank is set to \(16\), while GEAR initializes at \(r_0=4\), with a weighted-average rank budget of \(16\) and a per-module upper bound \(r_{\max}=32\); the final ranks are determined by the accepted expansions rather than the upper bound. Performance is reported together with the number of trainable parameters to evaluate both adaptation quality and parameter efficiency.

As shown in Table~\ref{tab:commonsense_results}, GEAR achieves the highest average performance on both backbones with substantially fewer trainable parameters. On LLaMA3-8B, it reaches \(88.45\), outperforming FlexLoRA by \(1.84\) points with \(14.4\)M versus \(42.0\)M parameters; on Qwen2.5-7B, it attains \(89.40\) with \(14.0\)M parameters, versus \(87.65\) with \(40.3\)M for FlexLoRA. Unlike adaptive-rank methods that rely mainly on statistics of the current low-rank representation, GEAR expands the adaptation space according to task-relevant normal directions while optimizing accessible directions on the current manifold, yielding a better performance--parameter trade-off across both backbones.
\begin{table*}[!t]
\centering
\caption{
Commonsense reasoning performance on LLaMA3-8B and Qwen2.5-7B. We report the mean over three runs with different random seeds. Best results are shown in \textbf{bold}.
}
\label{tab:commonsense_results}

\setlength{\tabcolsep}{3.6pt}
\renewcommand{\arraystretch}{1}

\begin{tabular}{l cccccccc c}
\toprule
Method (Params)
& BoolQ
& PIQA
& SIQA
& HellaS.
& WinoG.
& ARC-e
& ARC-c
& OBQA
& Avg. \\
\midrule

\multicolumn{10}{l}{\textit{\textbf{LLaMA3-8B}}} \\
\addlinespace[1pt]

LoRA (41.9M)
& 72.37 & 87.83 & 81.34 & 94.81 & 85.08 & 91.46 & 81.32 & 86.00 & 85.03 \\

AdaLoRA (31.6M)
& 72.23 & 87.21 & 81.27 & 94.60 & 85.16 & 93.18 & 82.68 & 85.80 & 85.27 \\

Riemannion (41.9M)
& 72.45 & 87.32 & 81.17 & 94.95 & 87.64 & 93.14 & 82.51 & 87.40 & 85.82 \\

GeoLoRA (41.3M)
& 72.84 & 87.54 & 79.99 & 94.62 & 85.56 & 93.52 & 83.96 & 87.00 & 85.63 \\

FlexLoRA (42.0M)
& 72.91 & 88.30 & 82.96 & 95.68 & 87.14 & 93.56 & 83.70 & 88.60 & 86.61 \\

GEAR (\textbf{14.4M})
& \textbf{75.19}
& \textbf{90.81}
& \textbf{83.78}
& \textbf{96.82}
& \textbf{89.82}
& \textbf{93.77}
& \textbf{85.58}
& \textbf{91.80}
& \textbf{88.45} \\

\midrule

\multicolumn{10}{l}{\textit{\textbf{Qwen2.5-7B}}} \\
\addlinespace[1pt]

LoRA (40.4M)
& 71.54 & 88.52 & 78.73 & 94.05 & 82.56 & 96.13 & 90.00 & 90.12 & 86.46 \\

AdaLoRA (35.2M)
& 71.59 & 89.12 & 78.92 & 94.06 & 82.79 & 96.00 & 90.27 & 89.80 & 86.57 \\

Riemannion (40.4M)
& 71.68 & 89.17 & 79.32 & 94.72 & 84.53 & 95.96 & \textbf{90.53} & 91.00 & 87.11 \\

GeoLoRA (40.1M)
& 71.90 & 89.12 & 79.02 & 94.88 & 85.00 & 96.21 & 89.93 & 91.00 & 87.13 \\

FlexLoRA (40.3M)
& 72.29 & 89.17 & 79.68 & 95.19 & 86.03 & 96.21 & 90.19 & 92.40 & 87.65 \\

GEAR (\textbf{14.0M})
& \textbf{75.51}
& \textbf{90.86}
& \textbf{81.58}
& \textbf{96.34}
& \textbf{90.53}
& \textbf{96.80}
& 89.76
& \textbf{93.80}
& \textbf{89.40} \\

\bottomrule
\end{tabular}
\end{table*}
\subsection{Natural Language Understanding}
We evaluate GEAR with DeBERTa-v3-base \citep{he2023debertav3} on the GLUE benchmark \citep{wang2019glue}, covering eight tasks with linguistic objectives, including acceptability, inference, paraphrase detection, semantic similarity, and sentiment classification. The reference rank is set to \(8\), while GEAR initializes at \(r_0=4\) with \(r_{\max}=16\) as an upper bound; the final rank is determined by the accepted expansions rather than the limit. We report the task-specific metrics and their average.

As shown in Table~\ref{tab:glue_results}, GEAR achieves the best average score of \(89.63\) with only \(1.5\)M trainable parameters, improving over FlexLoRA by \(0.78\) points while using about \(21\%\) fewer parameters. GEAR attains the strongest result on seven of the eight tasks, with particularly clear gains on MRPC, RTE, and STS-B. The improvement across heterogeneous objectives suggests that the benefit of GEAR is not tied to a specific task structure. By expanding the adaptation space according to task-relevant directions beyond its current first-order reach, GEAR can introduce additional capacity more selectively across diverse tasks.

\begin{table*}[!t]
\centering
\caption{
GLUE performance with DeBERTa-v3-base. We report the mean over three runs with different random seeds. Best results are shown in \textbf{bold}.
}
\label{tab:glue_results}

\setlength{\tabcolsep}{4.9pt}
\renewcommand{\arraystretch}{1}

\begin{tabular}{l cccccccc c}
\toprule
Method (Params)
& \shortstack{CoLA\\Mcc.}
& \shortstack{MNLI\\Acc.}
& \shortstack{MRPC\\Acc.}
& \shortstack{RTE\\Acc.}
& \shortstack{QNLI\\Acc.}
& \shortstack{SST-2\\Acc.}
& \shortstack{STS-B\\Corr.}
& \shortstack{QQP\\Acc.}
& Avg. \\
\midrule

LoRA (1.9M)
& 68.31 & 89.74 & 90.19 & 84.16
& 93.61 & 95.30 & 85.20 & 90.35 & 87.11 \\

AdaLoRA (1.9M)
& 70.07 & 89.93 & 90.44 & 88.17
& \textbf{94.27} & 95.13 & 91.21 & \textbf{91.41} & 88.83 \\

Riemannion (1.9M)
& 67.41 & 88.93 & 89.01 & 85.13
& 92.64 & 94.00 & 90.12 & 87.92 & 86.90 \\

GeoLoRA (1.9M)
& 68.18 & 89.67 & 87.76 & 87.73
& 92.93 & 95.52 & 90.15 & 87.77 & 87.46 \\

FlexLoRA (1.9M)
& 70.22 & 90.13 & 90.67 & 88.45
& 94.24 & 95.34 & 91.53 & 90.27 & 88.85 \\

GEAR (\textbf{1.5M})
& \textbf{70.85}
& \textbf{90.97}
& \textbf{92.46}
& \textbf{89.53}
& \textbf{94.27}
& \textbf{96.22}
& \textbf{91.83}
& 90.88
& \textbf{89.63} \\

\bottomrule
\end{tabular}
\end{table*}

\subsection{Visual Task}
To assess the generality of GEAR beyond NLP, we evaluate it with ViT-B/16 \citep{dosovitskiy2021image} on the Visual Task Adaptation Benchmark (VTAB-1K) \citep{zhai2020vtab}, which covers 19 recognition tasks grouped into Natural, Specialized, and Structured domains. We use the same rank configuration as in the NLU experiments.

As shown in Table~\ref{tab:vtab_results}, GEAR achieves the highest overall average of \(73.48\) with only \(0.2\)M trainable parameters, compared with \(72.55\) for FlexLoRA using \(0.3\)M. GEAR obtains the best result on 14 of the 19 tasks, including 7 of the 8 Structured tasks. The gains across visually and semantically diverse domains indicate that the proposed rank-evolution mechanism is not tied to a particular modality or task family. By adapting the low-rank space according to task-relevant directions beyond its current first-order reach, GEAR maintains strong transfer performance while introducing fewer additional degrees of freedom.
\begin{table}[t]
\centering
\caption{
VTAB-1K performance with ViT-B/16. We report the mean over three runs with different random seeds.
Best results are shown in \textbf{bold}.
}
\label{tab:vtab_results}

\setlength{\tabcolsep}{2.0pt}
\renewcommand{\arraystretch}{0.98}

\begin{tabular*}{\columnwidth}{@{\extracolsep{\fill}}lcccccc@{}}
\toprule
\multirow{2}{*}{\textbf{Dataset}}
& \textbf{LoRA}
& \textbf{AdaLoRA}
& \textbf{Riemannion}
& \textbf{GeoLoRA}
& \textbf{FlexLoRA}
& \textbf{GEAR} \\
&
{(0.3M)}
&
{(0.3M)}
&
{(0.3M)}
&
{(0.3M)}
&
{(0.3M)}
&
{(\textbf{0.2M})} \\
\midrule

\multicolumn{7}{l}{\textit{Natural}} \\
CIFAR100  & 69.39 & 71.81 & 70.78 & 68.05 & 71.30 & \textbf{71.88} \\
Caltech101 & 90.60 & 89.96 & 91.03 & 88.87 & 91.72 & \textbf{91.78} \\
DTD        & 71.01 & \textbf{71.76} & 71.28 & 68.88 & 70.18 & 71.60 \\
Flowers    & 98.80 & 99.04 & 99.09 & 99.06 & 98.98 & \textbf{99.12} \\
Pets       & 89.89 & 91.58 & 91.85 & 91.71 & 91.14 & \textbf{91.93} \\
SVHN       & 90.53 & 90.59 & 89.08 & 91.39 & 90.25 & \textbf{91.89} \\
SUN397     & 53.78 & \textbf{55.61} & 55.12 & 54.63 & 54.14 & 55.03 \\
\midrule

\multicolumn{7}{l}{\textit{Specialized}} \\
Camelyon    & 83.26 & 81.92 & 81.60 & 82.16 & \textbf{83.98} & 83.93 \\
EuroSAT     & 93.74 & 94.04 & 93.13 & 94.61 & 93.91 & \textbf{94.91} \\
RESISC45    & 85.40 & 85.37 & 84.30 & 84.52 & 85.41 & \textbf{85.98} \\
Retinopathy & 73.60 & 75.40 & 75.16 & \textbf{75.41} & 74.13 & 74.91 \\
\midrule

\multicolumn{7}{l}{\textit{Structured}} \\
CLEVR-C    & 80.91 & 81.23 & 78.95 & 80.63 & 81.30 & \textbf{82.07} \\
CLEVR-D    & 64.74 & 65.38 & 68.49 & 63.37 & 63.56 & \textbf{68.51} \\
DMLab      & 50.04 & 50.23 & 49.95 & 50.08 & 52.02 & \textbf{52.18} \\
KITTI-D    & 77.87 & 78.34 & 75.81 & 78.14 & \textbf{78.95} & 78.90 \\
dSpr-Loc   & 75.76 & 73.38 & 72.40 & 73.62 & 75.49 & \textbf{76.78} \\
dSpr-Ori   & 47.65 & 49.32 & 48.71 & 48.27 & 49.02 & \textbf{49.52} \\
sNORB-Azi  & 29.49 & 28.00 & 27.76 & 30.25 & 29.84 & \textbf{31.32} \\
sNORB-Ele  & 41.93 & 37.56 & 39.58 & 40.99 & 43.17 & \textbf{43.89} \\
\midrule

\textbf{Average}
& 72.02 & 72.13 & 71.79 & 71.82 & 72.55 & \textbf{73.48} \\
\bottomrule
\end{tabular*}
\end{table}

Across all evaluated settings, GEAR achieves the strongest average performance with fewer trainable parameters than the strongest adaptive-rank baseline. These gains span two language-model backbones, eight GLUE tasks, and 19 VTAB-1K tasks, indicating robustness across backbones, task families, and modalities. Importantly, they arise from adaptively selected final ranks rather than simply increasing the rank budget, supporting the hypothesis that identifying valuable directions outside the current adaptation space enables more efficient rank growth.

\section{Conclusion}
We presented GEAR, a geometry-guided framework for dynamic low-rank adaptation that shifts the focus from capacity allocation to the sufficiency of the current adaptation space. We established a first-order information boundary for standard factor gradients, showing that task-relevant normal information is unobservable from the current fixed-rank parameterization, and further related the normal spectrum to the locally optimal value of rank expansion. Building on this analysis, GEAR combines fixed-rank manifold optimization, efficient normal-spectrum estimation, and task-induced geometric initialization to adapt rank according to task-relevant directions beyond the current first-order reachable space. Across commonsense reasoning, GLUE, and VTAB-1K, GEAR achieves the highest average performance in all evaluated settings while using substantially fewer trainable parameters than strong adaptive-rank baselines. Overall, our results suggest that dynamic rank is most effectively viewed not as a budget-allocation problem alone, but as the controlled evolution of the adaptation space in response to previously unrepresented task information.

\section*{AI use statement}
Generative AI tools were used to assist with language editing, \LaTeX{}
formatting, and the examination of mathematical derivations. All AI-assisted
suggestions were independently reviewed and verified by the authors, who take
full responsibility for the final content of this work.

\section*{Reproducibility Statement}
We provide complete proofs of all theoretical results and detailed algorithm descriptions in the appendix. Experimental configurations, hyperparameters, dataset protocols, baseline settings, and additional ablations are also provided to facilitate reproducibility. All reported results are averaged over three random seeds.

\bibliography{bibfile}
\bibliographystyle{iclr2027_conference}



\newpage
\appendix
\numberwithin{equation}{subsection}
\renewcommand{\theequation}
{\Alph{section}.\arabic{subsection}.\arabic{equation}}


\section{Supplementary Overview}
\label{app:overview}

This supplementary material provides additional methodological, theoretical, and experimental details for GEAR. 
Appendix~\ref{app:method_details} describes the implementation of fixed-rank retraction and the validation procedure for rank-expansion candidates. 
Appendix~\ref{app:proofs} provides complete proofs of the three theoretical results presented in the main paper. 
Appendix~\ref{app:exp_details} reports the datasets, evaluation protocols, training configurations, and hyperparameter settings used in our experiments. Finally, ablation studies are provided in 
Appendix~\ref{app:ablation}.

\section{Additional Method Details}
\label{app:method_details}
\subsection{Low-Rank Retraction}
\label{app:detail_retraction}

We detail the retraction used in Section~\ref{sec:3.2}. Recall that the current adaptation state is represented as \(X=UC\), with \(U^\top U=I_r\), and let
\[
D = U^\top G,
\qquad
K=(I-UU^\top)GV .
\]
Using the tangent-gradient decomposition in Eq.~\eqref{eq:tangent_grad}, the update before retraction can be written as
\[
X^{-}
=
U(C-\eta D)-\eta K V^\top
=
LR,
\]
where
\[
L=[\,U,\;K\,],
\qquad
R=
\begin{bmatrix}
C-\eta D\\
-\eta V^\top
\end{bmatrix}.
\]
Hence, \(X^{-}\) has rank at most \(2r\) and can be manipulated without forming the full \(n\times m\) matrix.

We compute thin QR factorizations of the two compact factors,
\[
L=Q_LR_L,
\qquad
R^\top=Q_RR_R,
\]
which gives
\[
X^{-}
=
Q_L
\left(R_LR_R^\top\right)
Q_R^\top .
\]
The matrix \(R_LR_R^\top\) is at most \(2r\times2r\). We therefore compute only its compact SVD,
\[
R_LR_R^\top
=
P\Sigma Q^\top,
\]
and retain the leading \(r\) singular components. The retracted state is
\[
X_{t+1}
=
(Q_LP_r)\Sigma_r(Q_RQ_r)^\top .
\]
Accordingly, the updated low-rank representation can be obtained directly as
\[
U_{t+1}=Q_LP_r,
\qquad
V_{t+1}=Q_RQ_r,
\qquad
C_{t+1}=\Sigma_rV_{t+1}^\top .
\]
Thus, the retraction simultaneously restores the rank-\(r\) constraint and updates the left and right subspace bases required by the next iteration, while avoiding any SVD of the full adaptation matrix.

\subsection{Rank-Expansion Candidate Validation}
\label{app:detail_validate}

At each rank audit, the selected normal directions are treated as candidate
expansions rather than committed immediately. We order the proposed directions
by the utility score in Eq.~\eqref{eq:rank_score} and construct several
prefixes with different birth-step scales. The current adapter and all
candidate variants are evaluated on the same set of separate training microbatches under matched stochasticity, so that their losses are directly
comparable.

Let $\bar L_0$ and $\bar L_{\mathrm{cand}}$ denote the mean training losses of
the current adapter and a candidate, respectively. We define the relative
improvement as
\[
\mathrm{RI}
=
\frac{\bar L_0-\bar L_{\mathrm{cand}}}{|\bar L_0|}.
\]
A candidate is accepted only if $\mathrm{RI}>\rho$, where $\rho$ is the
prescribed acceptance threshold.

Among the accepted candidates, we favor the smallest rank increase whose loss
is close to the best observed value, and then the smaller birth-step scale.
The selected expansion is committed only if it also satisfies the per-module
rank limit and the global budget. If no candidate passes the acceptance
criterion, the audit leaves the adapter unchanged.
\subsection{Computational and Memory Complexity}
\label{app:complexity}
GEAR preserves the low-rank computational structure throughout training. 
For a rank-$r$ adapter, a standard LoRA factor update requires 
$\mathcal{O}((n+m)r)$ element-wise operations, while the additional retraction in GEAR costs
$\mathcal{O}((n+m)r^2+r^3)$ due to thin QR factorizations and an SVD of a $2r\times2r$ core matrix.
Since $r\ll \min\{n,m\}$, this overhead remains confined to the low-rank factors and small core matrices, without introducing any full $n\times m$ decomposition. 
The committed adapter state requires only $\mathcal{O}((n+m)r)$ memory.

TGI is performed once before fine-tuning. In the language experiments,
its $2(p+1)$ probing stages require $2(p+1)$ calibration backward
passes per initialization run and
$\mathcal O(h\max\{n,m\})$ peak sketch storage per adapted matrix.
The QR and SVD operations cost
$\mathcal O((p+1)(n+m)h^2)$.
We use $p=1$ and $h\ll\min\{n,m\}$. In VTAB, initialization instead aggregates dense gradient proxies for the selected Q/V projections from three calibration batches, requiring $\mathcal O(nm)$ temporary gradient storage per adapted matrix.

Rank audits introduce additional computation only at scheduled audit steps. 
With the same probe width $h$, the projections, orthogonalization, and decomposition of the compressed normal sketch incur
$\mathcal{O}((n+m)rh+(n+m)h^2)$ low-rank linear-algebra cost and 
$\mathcal{O}((n+m)h)$ temporary memory, in addition to the backward probes used to obtain gradient products. 
Both $r$ and $h$ are small relative to the ambient dimensions, and no full gradient or normal matrix is explicitly stored.

The main audit overhead therefore comes from the additional backward probes and candidate forward evaluations rather than large-scale matrix operations. 
These computations are performed only for $t\in\mathcal{T}_{\mathrm{audit}}$, while ordinary training steps incur only the fixed-rank manifold update. 
Moreover, all probing variables are transient and do not contribute to the committed adapter rank or the trainable-parameter counts reported in our experiments.

We further assess the measured resource costs on the DTD task with
ViT-B/16 in Table~\ref{tab:dtd_resource}.

\begin{table}[t]
\centering
\caption{Training time and peak CUDA memory on DTD with ViT-B/16.}
\label{tab:dtd_resource}
\begin{tabular}{lcc}
\toprule
Method & Training time (s) & Peak CUDA memory (MiB) \\
\midrule
LoRA     & 789.83  & 3825.15 \\
FlexLoRA & 1006.75 & 6553.87 \\
GEAR     & 1060.83 & 4458.21 \\
\bottomrule
\end{tabular}
\end{table}

\section{Proofs}
\label{app:proofs}

\subsection{Proof of Theorem \ref{thm:1}: Factor-Gradient Observability Barrier}
\label{app:proof_theorem1}

\begin{proof}
Consider the factorization map
\[
\Phi:\mathbb{R}^{n\times r}\times\mathbb{R}^{r\times m}
\rightarrow
\mathbb{R}^{n\times m},
\qquad
\Phi(B,A)=BA .
\]
Its differential at $(B,A)$ is
\[
D\Phi_{(B,A)}[\Delta B,\Delta A]
=
\Delta BA+B\Delta A .
\]

For arbitrary perturbations $(\Delta B,\Delta A)$, define
$B(t)=B+t\Delta B$ and $A(t)=A+t\Delta A$. Since the set of full-rank
matrices is open, both factors remain full rank for sufficiently small $t$.
Hence $X(t)=B(t)A(t)\in\mathcal M_r$ locally around $t=0$, and
\[
\dot X(0)
=
\Delta BA+B\Delta A
=
D\Phi_{(B,A)}[\Delta B,\Delta A].
\]
Therefore,
\[
\operatorname{Im}D\Phi_{(B,A)}
\subseteq
T_X\mathcal M_r.
\]

To determine the dimension of this image, consider the kernel condition
$\Delta BA+B\Delta A=0$. Since $A$ has full row rank, it admits a right
inverse $A_R$ satisfying $AA_R=I_r$. Right-multiplication by $A_R$ gives
$\Delta B=-B\Delta A A_R$, so $\Delta B=BK$ for some
$K\in\mathbb R^{r\times r}$. Substitution into the kernel condition yields
$B(KA+\Delta A)=0$. Since $B$ has full column rank,
$\Delta A=-KA$. Conversely, every pair $(BK,-KA)$ satisfies
$D\Phi_{(B,A)}[BK,-KA]=0$. Thus,
\[
\ker D\Phi_{(B,A)}
=
\{(BK,-KA):K\in\mathbb R^{r\times r}\}.
\]
The parametrization by $K$ is injective because $B$ has full column rank;
hence $\dim\ker D\Phi_{(B,A)}=r^2$. By the rank--nullity theorem,
\[
\dim\operatorname{Im}D\Phi_{(B,A)}
=
nr+rm-r^2
=
r(n+m-r).
\]
Since $\dim T_X\mathcal M_r=\dim\mathcal M_r=r(n+m-r)$, the inclusion above
is between spaces of equal dimension, and therefore
\[
\operatorname{Im}D\Phi_{(B,A)}
=
T_X\mathcal M_r.
\]

We next identify the adjoint of $D\Phi_{(B,A)}$. Under the Frobenius inner
product,
\[
\begin{aligned}
\left\langle
G,D\Phi_{(B,A)}[\Delta B,\Delta A]
\right\rangle_F
&=
\langle G,\Delta BA+B\Delta A\rangle_F \\
&=
\langle GA^\top,\Delta B\rangle_F
+
\langle B^\top G,\Delta A\rangle_F .
\end{aligned}
\]
Hence
\[
D\Phi_{(B,A)}^*(G)
=
(GA^\top,B^\top G)
=
\mathcal O_{A,B}(G).
\]

For any linear operator $L$ between finite-dimensional inner-product spaces,
$\ker L^*=(\operatorname{Im}L)^\perp$. Applying this identity to
$D\Phi_{(B,A)}$ and using the characterization above gives
\[
\ker\mathcal O_{A,B}
=
\left(\operatorname{Im}D\Phi_{(B,A)}\right)^\perp
=
\left(T_X\mathcal M_r\right)^\perp
=
N_X\mathcal M_r.
\]

Finally, for $D=G_1-G_2$,
$\mathcal O_{A,B}(G_1)=\mathcal O_{A,B}(G_2)$ if and only if
$\mathcal O_{A,B}(D)=0$, which is equivalent to $D\in N_X\mathcal M_r$.
Equivalently, the tangent projection of $D$ vanishes, and hence
\[
(G_1)_T=(G_2)_T.
\]
This proves both statements.
\end{proof}


\subsection{Proof of Theorem \ref{thm:2}: Local Optimality of Rank Expansion}
\label{app:proof_theorem2}

\begin{proof}
Let $X=U\Sigma V^\top\in\mathcal M_r$ be a compact SVD, and let
$U_\perp$ and $V_\perp$ be orthonormal complements of $U$ and $V$.
For any perturbation $H$, write
\[
\begin{bmatrix}
U^\top\\ U_\perp^\top
\end{bmatrix}
H
\begin{bmatrix}
V & V_\perp
\end{bmatrix}
=
\begin{bmatrix}
H_{11} & H_{12}\\
H_{21} & H_{22}
\end{bmatrix}.
\]
Its normal component is $H_N=U_\perp H_{22}V_\perp^\top$, while the remaining
three blocks constitute the tangent component $H_T$.

We first characterize the tangent cone of
$\mathcal V_{r+k}=\{Y:\operatorname{rank}(Y)\le r+k\}$ at $X$. In the above
coordinates,
\[
\mathcal C_{r,k}(X)
=
\left\{
H_T+H_N:
H_T\in T_X\mathcal M_r,\;
H_N\in N_X\mathcal M_r,\;
\operatorname{rank}(H_N)\le k
\right\}.
\]
To verify this characterization, consider a feasible curve
$X(t)\in\mathcal V_{r+k}$ with $X(0)=X$ and $\dot X(0)=H$. Its block form is
\[
X(t)
=
\begin{bmatrix}
\Sigma+tH_{11}+o(t) & tH_{12}+o(t)\\
tH_{21}+o(t) & tH_{22}+o(t)
\end{bmatrix}.
\]
For sufficiently small $t$, the upper-left block is nonsingular. Since
$\operatorname{rank}(X(t))\le r+k$, its Schur complement has rank at most
$k$. Dividing by $t$ and letting $t\to0$ therefore gives
$\operatorname{rank}(H_{22})\le k$.

Conversely, suppose $\operatorname{rank}(H_{22})\le k$ and define
$A_t=\Sigma+tH_{11}$. Consider
\[
\widetilde X(t)
=
\begin{bmatrix}
A_t & tH_{12}\\
tH_{21} &
tH_{22}+t^2H_{21}A_t^{-1}H_{12}
\end{bmatrix}.
\]
Its Schur complement with respect to $A_t$ is $tH_{22}$, so
$\operatorname{rank}(\widetilde X(t))
=r+\operatorname{rank}(H_{22})\le r+k$. Moreover,
$\widetilde X(0)=X$ and $\dot{\widetilde X}(0)=H$, establishing the stated
tangent-cone characterization.

We now minimize
\[
q_\eta(H)
=
\langle G,H\rangle_F
+
\frac{1}{2\eta}\|H\|_F^2
\]
over $\mathcal C_{r,k}(X)$. Since the tangent and normal spaces are orthogonal,
write $G=G_T+G_N$ and $H=H_T+H_N$. Then
\[
q_\eta(H)
=
\left(
\langle G_T,H_T\rangle_F
+
\frac{1}{2\eta}\|H_T\|_F^2
\right)
+
\left(
\langle G_N,H_N\rangle_F
+
\frac{1}{2\eta}\|H_N\|_F^2
\right),
\]
subject only to $\operatorname{rank}(H_N)\le k$.

For the tangent component,
\[
\langle G_T,H_T\rangle_F
+
\frac{1}{2\eta}\|H_T\|_F^2
=
\frac{1}{2\eta}\|H_T+\eta G_T\|_F^2
-
\frac{\eta}{2}\|G_T\|_F^2,
\]
so the unique minimizer is $H_T^\star=-\eta G_T$.

Similarly, for the normal component,
\[
\langle G_N,H_N\rangle_F
+
\frac{1}{2\eta}\|H_N\|_F^2
=
\frac{1}{2\eta}\|H_N+\eta G_N\|_F^2
-
\frac{\eta}{2}\|G_N\|_F^2.
\]
Thus, under $\operatorname{rank}(H_N)\le k$, minimizing the normal term is
equivalent to finding the best rank-$k$ approximation of $-\eta G_N$ in
Frobenius norm. By the Eckart--Young--Mirsky theorem,
$H_N^\star=-\eta(G_N)_k$. Hence
\[
H_k^\star
=
-\eta\left[G_T+(G_N)_k\right].
\]

It remains to evaluate the additional decrease relative to fixed-rank
optimization. Over $T_X\mathcal M_r$,
$\min q_\eta=-\frac{\eta}{2}\|G_T\|_F^2$, whereas over
$\mathcal C_{r,k}(X)$,
\[
\min q_\eta
=
-\frac{\eta}{2}\|G_T\|_F^2
-
\frac{\eta}{2}\|(G_N)_k\|_F^2.
\]
Therefore,
\[
\Delta_k
=
\frac{\eta}{2}\|(G_N)_k\|_F^2
=
\frac{\eta}{2}\sum_{i=1}^{k}\sigma_i^2(G_N).
\]
Thus, the leading normal singular directions yield the optimal local
rank-$k$ expansion, and their squared singular values determine the additional
decrease of the local quadratic model.
\end{proof}
\subsection{Proof of Theorem \ref{thm:3}: Stability under Inexact Normal Estimation}
\label{app:proof_theorem3}

\begin{proof}
For brevity, let $N=G_N$, $\widehat N=\widehat G_N$,
$N_k=(G_N)_k$, and $\widehat N_k=(\widehat G_N)_k$.
By assumption, $\|\widehat N-N\|_F\le\varepsilon$ and
$\tau_k=\|N-N_k\|_F$.

Since $N_k$ and $\widehat N_k$ lie in the normal space, they are orthogonal
to the tangent component of $G$. Hence, for any normal matrix $B$,
\[
q_\eta(-\eta B)
=
\frac{\eta}{2}\|N-B\|_F^2
-
\frac{\eta}{2}\|N\|_F^2.
\]
Applying this identity to $B=N_k$ and $B=\widehat N_k$ yields
\[
q_\eta(-\eta\widehat N_k)-q_\eta(-\eta N_k)
=
\frac{\eta}{2}
\left(
\|N-\widehat N_k\|_F^2-\|N-N_k\|_F^2
\right).
\]
Because $N_k$ is a best rank-$k$ approximation of $N$ and
$\operatorname{rank}(\widehat N_k)\le k$,
$\|N-N_k\|_F\le\|N-\widehat N_k\|_F$. Therefore,
$q_\eta(-\eta\widehat N_k)-q_\eta(-\eta N_k)\ge0$.

It remains to establish the upper bound. Since $\widehat N_k$ is a best
rank-$k$ approximation of $\widehat N$,
$\|\widehat N-\widehat N_k\|_F\le\|\widehat N-N_k\|_F$. Using the triangle
inequality,
\[
\begin{aligned}
\|N-\widehat N_k\|_F
&\le
\|N-\widehat N\|_F
+
\|\widehat N-\widehat N_k\|_F \\
&\le
\varepsilon+\|\widehat N-N_k\|_F \\
&\le
\varepsilon+\|\widehat N-N\|_F+\|N-N_k\|_F \\
&\le
\tau_k+2\varepsilon .
\end{aligned}
\]
Substituting this bound into the preceding identity gives
\[
\begin{aligned}
q_\eta(-\eta\widehat N_k)-q_\eta(-\eta N_k)
&\le
\frac{\eta}{2}
\left[(\tau_k+2\varepsilon)^2-\tau_k^2\right] \\
&=
2\eta\varepsilon(\tau_k+\varepsilon).
\end{aligned}
\]
Consequently,
\[
0
\le
q_\eta\!\left(-\eta(\widehat G_N)_k\right)
-
q_\eta\!\left(-\eta(G_N)_k\right)
\le
2\eta\varepsilon(\tau_k+\varepsilon).
\]
This completes the proof.
\end{proof}

\section{Experimental Details}
\label{app:exp_details}

\subsection{Details on Commonsense Reasoning Task}
\label{app:exp_commonsense}
The commonsense reasoning suite comprises eight benchmarks: BoolQ
\citep{clark2019boolq}, PIQA \citep{bisk2020piqa}, SIQA
\citep{sap2019socialiqa}, HellaSwag \citep{zellers2019hellaswag}, WinoGrande
\citep{sakaguchi2020winogrande}, ARC-e/ARC-c \citep{clark2018arc}, and OBQA
\citep{mihaylov2018openbookqa}. Following the standard multi-task protocol, we
combine the training portions of these tasks into the Commonsense170K corpus \citep{hu2023llmadapters}
and evaluate each benchmark separately on its corresponding test set. We use
LLaMA3-8B and Qwen2.5-7B as backbone models. The main hyperparameter settings
are summarized in Table~\ref{tab:commonsense_hyperparameters}.

The probing configuration estimates task-induced gradient structure for geometric initialization. We use probe width $h=24$ for both TGI and rank audits.
For rank audits, each calibration batch contains 128 examples,
with 16 from each task. The damping coefficient $\lambda=0.1$
regularizes the gradient-scale normalization. Subsequent experiments
follow the same initialization and rank-audit framework with
task-specific adjustments; the VTAB implementation is detailed in
Appendix~\ref{app:exp_visual}.

\begin{table}[htpb]
\centering
\caption{
Hyperparameter settings for LLaMA3-8B and Qwen2.5-7B in commonsense reasoning experiments.
}
\label{tab:commonsense_hyperparameters}

\setlength{\tabcolsep}{6pt}
\renewcommand{\arraystretch}{1.05}

\begin{tabular}{@{}ll@{}}
\toprule
\textbf{Hyperparameter} & \textbf{Setting} \\
\midrule

Initial effective rank $r_0$        & 4 \\
Maximum per-module rank $r_{\max}$  & 32 \\
Weighted-average rank budget $\bar r_{\mathrm{budget}}$       & 16 \\
Target modules                      & Q, K, V, O, Up, Down, Gate \\
\midrule

Probe $\Omega$ width $h$ & 24 \\
Probe batch size & 128 total (16 per task) \\
Damping $\lambda$                   & 0.1 \\
\midrule

First audit step $t_{\mathrm{start}}$  & 100 \\
Audit interval $\Delta t$              & 50 \\
Rank-audit cutoff $t_{\mathrm{final}}$ & 500 \\
Max. rank increment per module per audit                       & 4 \\
Acceptance threshold $\rho$            & 0.5\% \\
\midrule

Learning rate                       & $2\times10^{-2}$ \\
Dropout                             & 0.05 \\
Scheduler / warmup                  & Linear / 0.1 \\
Batch configuration                 & 16 microbatch $\times$ 4 accumulation \\
Epochs                              & 2 \\

\bottomrule
\end{tabular}
\end{table}

\subsection{Details on Natural Language Understanding Task}
\label{app:exp_nlu}
For the natural language understanding experiments, we fine-tune DeBERTa-v3-base on the GLUE benchmark \citep{wang2019glue}, comprising eight tasks across
single-sentence classification, semantic similarity and paraphrase detection,
and natural language inference. Specifically, we evaluate on CoLA
\citep{warstadt2019cola}, SST-2 \citep{socher2013sst}, MRPC
\citep{dolan2005mrpc}, QQP \citep{wang2019glue}, STS-B
\citep{cer2017stsb}, MNLI \citep{williams2018mnli}, QNLI
\citep{rajpurkar2016squad}, and RTE
\citep{dagan2006rte,barhaim2006rte,giampiccolo2007rte,bentivogli2009rte}.
Dataset statistics and evaluation metrics are summarized in
Table~\ref{tab:glue_dataset}, and the corresponding training and
rank-adaptation settings are reported in
Table~\ref{tab:glue_hyperparameters}. We use the same probing and subspace estimation procedure as in the commonsense reasoning experiments, with shared probe width $h=24$ for TGI and rank audits and task-specific calibration data and optimization settings. Separate learning rates are used for the adapter and task head because they operate under different parameterizations and effective update scales: the head follows standard Euclidean optimization, whereas the adapter learning rate controls the manifold-constrained update before retraction.

\begin{table*}[htpb]
\centering
\caption{
Details of the GLUE dataset.
}
\label{tab:glue_dataset}

\setlength{\tabcolsep}{5pt}
\renewcommand{\arraystretch}{1.10}

\begin{tabular}{l|lccccc}
\toprule
\textbf{Dataset}
& \textbf{Task}
& \textbf{\# Train}
& \textbf{\# Dev}
& \textbf{\# Test}
& \textbf{\# Label}
& \textbf{Metrics} \\
\midrule

\multicolumn{7}{c}{\textbf{Single-Sentence Classification}} \\
\midrule
CoLA
& Acceptability
& 8.5k
& 1k
& 1k
& 2
& Matthews corr. \\

SST-2
& Sentiment
& 67k
& 872
& 1.8k
& 2
& Accuracy \\

\midrule
\multicolumn{7}{c}{\textbf{Similarity and Paraphrase}} \\
\midrule
MRPC
& Paraphrase
& 3.7k
& 408
& 1.7k
& 2
& Accuracy \\

QQP
& Paraphrase
& 364k
& 40k
& 391k
& 2
& Accuracy \\

STS-B
& Similarity
& 7k
& 1.5k
& 1.4k
& 1
& Pearson/Spearman corr. \\

\midrule
\multicolumn{7}{c}{\textbf{Natural Language Inference}} \\
\midrule
MNLI
& NLI
& 393k
& 20k
& 20k
& 3
& Accuracy \\

QNLI
& QA/NLI
& 108k
& 5.7k
& 5.7k
& 2
& Accuracy \\

RTE
& NLI
& 2.5k
& 276
& 3k
& 2
& Accuracy \\

\bottomrule
\end{tabular}
\end{table*}

\begin{table*}[htpb] \centering \caption{Hyper-parameter settings of our method on the GLUE benchmark.} \label{tab:glue_hyperparameters} \resizebox{\textwidth}{!}{ \begin{tabular}{lcccccccccc} \toprule \textbf{Dataset} & \textbf{Adapter LR} & \textbf{Head LR} & \textbf{Batch} & \textbf{Epochs} & $\boldsymbol{r_0}$ & $\boldsymbol{\bar r_{\mathrm{budget}}}$ & $\boldsymbol{r_{\max}}$ & $\boldsymbol{t_{\mathrm{start}}}$ & $\boldsymbol{\Delta t}$ & $\boldsymbol{t_{\mathrm{final}}}$ \\ \midrule MNLI & $5\!\times\!10^{-4}$ & $5\!\times\!10^{-4}$ & 32 & 6 & 4 & 8 & 16 & $0.1T$ & $0.075T$ & $0.7T$ \\ RTE & $3.6\!\times\!10^{-3}$ & $1.2\!\times\!10^{-3}$ & 32 & 50 & 4 & 8 & 16 & $0.1T$ & $0.075T$ & $0.7T$ \\ QNLI & $3\!\times\!10^{-3}$ & $3\!\times\!10^{-4}$ & 32 & 5 & 4 & 8 & 16 & $0.1T$ & $0.075T$ & $0.7T$ \\ MRPC & $3\!\times\!10^{-3}$ & $1\!\times\!10^{-3}$ & 32 & 30 & 4 & 8 & 16 & $0.1T$ & $0.075T$ & $0.7T$ \\ QQP & $1.6\!\times\!10^{-3}$ & $8\!\times\!10^{-4}$ & 32 & 5 & 4 & 8 & 16 & $0.1T$ & $0.075T$ & $0.7T$ \\ SST-2 & $5\!\times\!10^{-3}$ & $3\!\times\!10^{-4}$ & 32 & 24 & 4 & 8 & 16 & $0.1T$ & $0.075T$ & $0.7T$ \\ CoLA & $1.6\!\times\!10^{-3}$ & $8\!\times\!10^{-4}$ & 32 & 20 & 4 & 8 & 16 & $0.1T$ & $0.075T$ & $0.7T$ \\ STS-B & $4.4\!\times\!10^{-3}$ & $2.2\!\times\!10^{-3}$ & 32 & 25 & 4 & 8 & 16 & $0.1T$ & $0.075T$ & $0.7T$ \\ \bottomrule \end{tabular} } \vspace{1mm} \begin{minipage}{0.98\textwidth} \footnotesize $T$ denotes the total number of optimization steps. At each rank audit, each module may add at most one new rank-one direction, with an acceptance threshold of $\rho=0.1\%$ relative loss improvement. \end{minipage} \end{table*}

\subsection{Details on Visual Task}
\label{app:exp_visual}

We evaluate GEAR with ViT-B/16 on VTAB-1K \citep{zhai2020vtab}, which contains 19 visual classification tasks across Natural, Specialized, and Structured categories. Dataset statistics and hyperparameter settings are reported in Tables~\ref{tab:vtab_dataset_full} and~\ref{tab:vtab_hyperparameters}, respectively. For TGI, we aggregate dense gradient proxies for the selected Q/V projections from three calibration batches. Rank audits and candidate validation use task-specific settings.

\begin{table}[htpb]
\centering
\caption{Details of the VTAB-1K benchmark.}
\label{tab:vtab_dataset_full}

\setlength{\tabcolsep}{7pt}
\renewcommand{\arraystretch}{1.05}

\begin{tabular}{lrrrr}
\toprule
\textbf{Dataset}
& \textbf{\# Classes}
& \textbf{\# Train}
& \textbf{\# Val}
& \textbf{\# Test} \\
\midrule

\multicolumn{5}{l}{\textit{Natural}} \\
CIFAR100        & 100 & 800/1000 & 200 & 10,000 \\
Caltech101      & 102 & 800/1000 & 200 & 6,084 \\
DTD             & 47  & 800/1000 & 200 & 1,880 \\
Flowers102      & 102 & 800/1000 & 200 & 6,149 \\
Pets            & 37  & 800/1000 & 200 & 3,669 \\
SVHN            & 10  & 800/1000 & 200 & 26,032 \\
SUN397          & 397 & 800/1000 & 200 & 21,750 \\

\midrule
\multicolumn{5}{l}{\textit{Specialized}} \\
Camelyon        & 2   & 800/1000 & 200 & 32,768 \\
EuroSAT         & 10  & 800/1000 & 200 & 5,400 \\
RESISC45        & 45  & 800/1000 & 200 & 6,300 \\
Retinopathy     & 5   & 800/1000 & 200 & 42,670 \\

\midrule
\multicolumn{5}{l}{\textit{Structured}} \\
CLEVR-Count     & 8   & 800/1000 & 200 & 15,000 \\
CLEVR-Dist      & 6   & 800/1000 & 200 & 15,000 \\
DMLab           & 6   & 800/1000 & 200 & 22,735 \\
KITTI-Dist      & 4   & 800/1000 & 200 & 711 \\
dSpr-Loc        & 16  & 800/1000 & 200 & 73,728 \\
dSpr-Ori        & 16  & 800/1000 & 200 & 73,728 \\
sNORB-Azim      & 18  & 800/1000 & 200 & 12,150 \\
sNORB-Ele       & 18  & 800/1000 & 200 & 12,150 \\

\bottomrule
\end{tabular}
\end{table}

\begin{table}[htpb]
\centering
\caption{
Hyperparameter settings for the VTAB-1K experiments.
}
\label{tab:vtab_hyperparameters}

\begin{tabular}{ll}
\toprule
\textbf{Hyperparameter} & \textbf{Setting} \\
\midrule
Initial effective rank $r_0$                  & 4 \\
Weighted-average rank budget $\bar r_{\mathrm{budget}}$                 & 8 \\
Maximum per-module rank $r_{\max}$            & 16 \\
Learning rate                                 & $5\times10^{-3}$ \\
First audit step $t_{\mathrm{start}}$          & $160$ \\
Audit interval $\Delta t$                     & $50$ \\
Rank-audit cutoff $t_{\mathrm{final}}$         & $1300$ \\
Acceptance threshold $\rho$                   & 0.03\% \\
Batch size                                    & 64 \\
Epochs                                        & 100 \\
\bottomrule
\end{tabular}
\end{table}

\section{Ablation Studies}
\label{app:ablation}
To isolate the contribution of each component, we conduct a component-wise
ablation on LLaMA3-8B over the commonsense reasoning suite. Starting from the
standard LoRA baseline, we examine task-induced geometric initialization
(TGI), fixed-rank manifold optimization, and geometry-guided rank expansion.
The results are reported in Table~\ref{tab:ablation_components}.

TGI and manifold optimization independently improve the average score from
$85.03$ to $85.57$ and $85.65$, respectively, while combining them further
raises it to $86.23$. This suggests that task-informed initialization and
geometry-preserving optimization provide complementary benefits in constructing
and optimizing the current adaptation space.

Rank expansion brings a larger gain. Adding it to manifold optimization improves the average from \(85.65\) to \(87.32\), while combining all three components further increases the score to \(88.45\), with improvements across all eight tasks. TGI also remains beneficial when rank expansion is enabled, improving the manifold-plus-expansion configuration from \(87.32\) to \(88.45\). Taken together, these ablations align with the central view of GEAR: dynamic rank adaptation should distinguish between better exploiting the current adaptation space and expanding it when task-relevant directions remain unrepresented.

\begin{table*}[htpb]
\centering
\caption{
Component ablation of GEAR on commonsense reasoning with LLaMA3-8B.
TGI denotes task-induced geometric initialization.
}
\label{tab:ablation_components}

\setlength{\tabcolsep}{2pt}
\renewcommand{\arraystretch}{1.05}

\begin{tabular}{ccc|ccccccccc}
\toprule
\textbf{TGI}
& \textbf{Manifold}
& \textbf{Rank Expansion}
& \textbf{BoolQ}
& \textbf{PIQA}
& \textbf{SIQA}
& \textbf{HellaS.}
& \textbf{WinoG.}
& \textbf{ARC-e}
& \textbf{ARC-c}
& \textbf{OBQA}
& \textbf{Avg.} \\
\midrule

-- & -- & --
& 72.37 & 87.83 & 81.34 & 94.81
& 85.08 & 91.46 & 81.32 & 86.00 & 85.03 \\

\checkmark & -- & --
& 72.51 & 87.57 & 81.27 & 94.88
& 86.94 & 92.16 & 83.08 & 86.16 & 85.57 \\

-- & \checkmark & --
& 72.40 & 87.03 & 81.12 & 94.46
& 87.38 & 92.51 & 83.12 & 87.21 & 85.65 \\

\checkmark & \checkmark & --
& 73.15 & 89.13 & 82.01 & 94.90
& 87.27 & 92.88 & 83.30 & 87.17 & 86.23 \\

-- & \checkmark & \checkmark
& 73.95 & 89.73 & 83.13 & 95.84
& 88.66 & 93.04 & 84.71 & 89.51 & 87.32 \\

\checkmark & \checkmark & \checkmark
& \textbf{75.19}
& \textbf{90.81}
& \textbf{83.78}
& \textbf{96.82}
& \textbf{89.82}
& \textbf{93.77}
& \textbf{85.58}
& \textbf{91.80}
& \textbf{88.45} \\

\bottomrule
\end{tabular}
\end{table*}


\end{document}